\documentclass{article}
\usepackage{amscd}
\usepackage{amsmath,amssymb,amsthm}
\usepackage{pifont,mathrsfs}
\usepackage{array,lastpage}
\usepackage{enumerate,xspace,pifont,shadow}
\usepackage{mathrsfs}
\usepackage{graphicx}
\usepackage{multimedia}
\usepackage{calrsfs}
\usepackage[all]{xy}

\DeclareSymbolFont{AMSb}{U}{msb}{m}{n}
\DeclareMathSymbol{\Z}{\mathbin}{AMSb}{"5A}
\DeclareMathSymbol{\R}{\mathbin}{AMSb}{"52}
\DeclareMathSymbol{\N}{\mathbin}{AMSb}{"4E}
\DeclareMathSymbol{\Q}{\mathbin}{AMSb}{"51}

\newcommand{\Th}{\textup{Th}}

\newcommand{\End}{\textup{End}}
\newcommand{\Aut}{\textup{Aut}}

\def\Ind{\setbox0=\hbox{$x$}\kern\wd0\hbox to 0pt{\hss$\mid$\hss}
\lower.9\ht0\hbox to 0pt{\hss$\smile$\hss}\kern\wd0}

\def\Notind{\setbox0=\hbox{$x$}\kern\wd0\hbox to 0pt{\mathchardef
\nn=12854\hss$\nn$\kern1.4\wd0\hss}\hbox to
0pt{\hss$\mid$\hss}\lower.9\ht0 \hbox to
0pt{\hss$\smile$\hss}\kern\wd0}

\newtheorem{sthm}{Theorem}
\newtheorem{thm}{Theorem}[section]
\newtheorem{lem}[thm]{Lemma}
\newtheorem{cor}[thm]{Corollary}

\newtheorem{fact}[thm]{Fact}

\theoremstyle{definition}
\newtheorem{definition}[thm]{Definition}

\theoremstyle{remark}

\theoremstyle{remark}

\theoremstyle{remark}
\newtheorem{claim}[thm]{Claim}

\begin{document}

\bibliographystyle{plain}

\title{The Schr\"{o}der-Bernstein Property for Theories of Abelian Groups}
\author{John Goodrick}
\maketitle

\begin{abstract}
A first-order theory has the \emph{Schr\"{o}der-Bernstein Property} if any two of its models that are elementarily bi-embeddable are isomorphic.  We prove (as Theorem~\ref{abeliangroups}):

\begin{sthm}
If $G$ is an abelian group, then the following are equivalent:

1. $\Th(G, +)$ has the Schr\"{o}der-Bernstein~property;

2. $\Th(G, +)$ is $\omega$-stable;

3. $G$ is the direct sum of a divisible group and a torsion group of bounded exponent;

4. $\Th(G, +)$ is superstable, and if $(\overline{G}, +) \equiv (G, +)$ is saturated, every map in $\Aut(\overline{G} / \overline{G}^\circ)$ is unipotent.
\end{sthm}

\end{abstract}

\section{Introduction}

A category $\textbf{C}$ has the \emph{Schr\"{o}der-Bernstein property}, or the ``SB~property'' for short, if any two of its objects that are bi-embeddable via monic morphisms are isomorphic.  (A morphism $f: A \rightarrow B$ is \emph{monic} if whenever $g_0, g_1 : C \rightarrow A$ and $f \circ g_0 = f \circ g_1$ then $g_0 = g_1$.)  If $T$ is a theory in first-order logic, then we say that $T$ has the SB~property if the category whose object class is $\textup{Mod}(T)$ and whose arrows are all elementary embeddings has the SB~property.  Note that in this category \emph{all} maps are injective, so the SB~property for $T$ is the same as any two elementarily bi-embeddable models being isomorphic.

The first investigation of which theories have the SB~property seems to have been done by Nurmagambetov in the papers \cite{nur2} and \cite{nur1}, in which he shows, among other things, that a countable $\omega$-stable theory has the SB~property if and only if it is nonmultidimensional.  In my thesis \cite{my_thesis}, I have shown that any countable first-order theory with the SB~property must be classifiable in the sense of Shelah (i.e. it must be superstable and have NOTOP and NDOP -- see also \cite{bible}) as well as nonmultidimensional.  The question of which superstable nonmultidimensional theories have the SB~property appears to be more delicate and no nice characterization of these theories has been proven; some conjectures about this are discussed in \cite{my_thesis} and \cite{when}.

The idea of focusing on the special case of abelian groups (suggested by my advisor Thomas Scanlon) comes from the hope that techniques from geometric stability theory may be used to understand the Schr\"{o}der-Bernstein property.  Very roughly, the idea is that checking whether a nonmultidimensional classifiable theory has the SB~property should reduce to examining the regular types, and the most interesting (and difficult) case seems to be when the regular type is nontrivial and locally~modular, in which case it is nonorthogonal to the generic type of a definable abelian group.

In this paper, we consider the question of which complete theories of abelian groups in the pure language of groups (i.e., $\left\{+\right\}$) have the SB~property.  Note that in this context elementary embeddings are the same thing as pure embeddings (see Lemma~\ref{purity}).  In the special case of such theories of abelian groups, SB~property turns out to be equivalent to $\omega$-stability (Theorem~\ref{abeliangroups} below).  Section~2 collects some relevant facts from the model~theory of modules and Section~3 contains the proof of the main theorem in the abstract.  Note that the direction 2 $\Rightarrow$ 1 follows already from Nurmagambetov's theorem, but we present here a new proof that $\omega$-stable abelian groups have the SB~property.  At a few points, we will assume familiarity with the standard model theoretic concepts of stability, superstability, and $\omega$-stability; for definitions of these terms, see \cite{baldwin}, \cite{hodges}, \cite{pillay}, or \cite{bible}.

\section{Theories of abelian groups}

This section contains some background material on theories of abelian groups in the pure group language.  Everything in this section has long been known by specialists in the area and we are just collecting a few results that we need for ease of reference.  For a comprehensive account of the subject, see the book \cite{prest}.

First, we set some notation.  Throughout this section, ``$p$'' will always denote a prime number in $\Z$.  $\Q$ is the group of all rational numbers and $\Z$ is the group of all integers, both with addition as the group operation.  For any prime $p$, ${\Z}_{p^{\infty}}$ denotes the subgroup of ${\Q}/{\Z}$ generated by the set $\left\{p^{-i} : i \in \omega \right\}$.  There are several other equivalent ways of describing this group: ${\Z}_{p^{\infty}}$ is isomorphic to the multiplicative group of all complex $p^n$-th roots of unity (for $p$ constant, $n$ varies), or it can be thought of as the direct limit of the diagram ${\Z}/{p\Z} \hookrightarrow {\Z}/{p^2\Z} \hookrightarrow {\Z}/{p^3\Z} \hookrightarrow \ldots$, where the embeddings are the unique nonzero maps between the groups.  

\begin{definition}
G is a \emph{$p$-group} if for every $g \in G$ there is an $n \in \omega$ such that $g$ has order $p^n$.
\end{definition}

So $\Z/{p \Z}$ and $\Z_{p^{\infty}}$ are both $p$-groups.

\begin{definition}
Let $G$ be an abelian group.

1. $G$ is \emph{divisible} if for any $g \in G$ and nonzero $n \in \omega$ there is a $y \in G$ such that $ny = x$.

2. $G$ is \emph{reduced} if $G$ contains no nonzero divisible subgroup.
\end{definition}

\begin{fact}
\label{div_class}
(\cite{kaplansky}, Theorem~4) Any divisible abelian group is a direct sum of copies of $\Q$ and $\Z_{p^{\infty}}$ for various primes $p$.  Moreover, the number of summands of each isomorphism type is a uniquely determined invariant of the group.
\end{fact}

\begin{fact}
\label{div-summand}
(Theorem~2 of \cite{kaplansky}) Any abelian group $G$ can be written as a direct sum $R \oplus D$ where $R$ is reduced and $D$ is divisible.
\end{fact}

\begin{definition}
1. If $G$ is an abelian group, then $H$ is a \emph{pure subgroup of $G$} if $H$ is a subgroup of $G$ and for any $h \in H$ and $n \in \omega$, if the equation $n x = h$ has a solution for $x$ in $G$ then it has a solution in $H$.

2. An injective group homomorphism $f: G \rightarrow H$ is a \emph{pure embedding} if its image $f(G)$ is a pure subgroup of $H$.
\end{definition}

Note that if $H$ is a direct summand of $G$ then $H$ is a pure subgroup of $G$.

From now on, we will freely use terminology and notation from model theory such as $\Th(G)$, elementary embeddings, saturated models, types, and ($\omega$-)stable theories.  Accordingly, if $(G, +)$ is an elementary submodel of $(H, +)$, we abbreviate this by $G \prec H$, and if $(G, +)$ is elementarily equivalent to $(H, +)$ then we write $G \equiv H$.

The most important facts about the model theory of abelian groups in the pure language of groups are the elimination of quantifiers up to p.p. formulas and the characterization of complete theories by Szmielew invariants. These are standard results that are explained in section~2.$\Z$ of \cite{prest}, and the next several facts are corollaries to these facts.

\begin{lem}
\label{purity}
If $f$ is an embedding between groups $G$ and $H$ that are both models of the same complete theory$T$, $f$ is pure if and only if it is elementary.
\end{lem}

\begin{proof}
Follows directly from Proposition~2.25 of \cite{prest}.
\end{proof}

\begin{definition}
An abelian group $G$ is \emph{pure-injective}, or \emph{algebraically~compact}, if every system of equations over $G$ (i.e. quantifier-free formulas) which is finitley satisfiable in $G$ has a solution in $G$.
\end{definition}

\begin{fact}
(Corollary~2.9 of \cite{prest}) Any $\omega_1$-saturated model of a theory $T$ as above is purely-injective, and so every model of $T$ is elementarily embeddable into a pure-injective group.
\end{fact}

We write $\Z_{(p)}$ for the localization of the ring $\Z$ at the ideal $(p)$.  The \emph{$p$-adic topology on $\Z_{(p)}$} is the topology whose basic open neighborhoods are the cosets of the subgroups $p^k \Z_{(p)}$ as $k$ varies.  In fact, this gives $\Z_{(p)}$ a metric structure such that each coset of $p^k \Z_{(p)}$ is a ball of radius $p^{-k}$, and the completion of $\Z_{(p)}$ as a metric space is denoted $\widehat{\Z_{(p)}}$ and called the set of $\emph{$p$-adic integers}$.  (We use the more cumbersome notation ``$\widehat{\Z_{(p)}}$'' rather than the standard ``$\Z_{(p)}$'' to avoid any possible confusion with the finite group $\Z / {p \Z}$, which is sometimes denoted $\Z_p$.)  When talking about groups, we will often abuse notation and use $\Z_{(p)}$ and $\widehat{\Z_{(p)}}$ to refer to the additive groups of the corresponding rings.  Sometimes we will be talking about elementary subgroups $G$ of $\widehat{\Z_{(p)}}$, in which case we write ``$\widehat{G}$'' for $\widehat{\Z_{(p)}}$.

\begin{thm}
\label{cases}
Any $T$ has a model that is a direct sum of direct-sum indecomposable (or simply ``indecomposable") pure-injective groups.  More concretely, it has a model that is a direct sum of finite cyclic groups, groups isomorphic to ${\Z}_{p^{\infty}}$, groups isomorphic to $\Q$, and groups isomorphic to $\widehat{{\Z}_{(p)}}$ (where $p$ may take on various values for different summands).
\end{thm}

For a proof and discussion of the last theorem, see Section~4.4 of \cite{prest} and Corollary~4.36 in particular.

\begin{fact}
\label{direct-sums}
Suppose that $\langle G_i : i \in I \rangle$ and $\langle H_i : i \in I \rangle$ are two sequences of abelian groups and for each $i \in I$, $G_i \equiv H_i$.  Then $\bigoplus_{i \in I} G_i \equiv \bigoplus_{i \in I} H_i$.
\end{fact}

\begin{proof}
This follows immediately from Lemma~2.23~(b) and (c) of \cite{prest}.
\end{proof}

\begin{thm}
\label{completions}
Suppose that $\langle G_i : i \in I \rangle$ is a set of elementary subgroups of $\widehat{\Z_{(p_i)}}$ for various primes $p_i$ and $H$ is any abelian group.  Then $H \oplus \bigoplus_{i \in I} G_i \prec H \oplus \bigoplus_{i \in I} \widehat{G_i}$.  (So in particular $H \oplus \bigoplus_{i \in I} G_i \equiv H \oplus \bigoplus_{i \in I} \widehat{G_i}$.)
\end{thm}

\begin{proof}
The fact that each $G_i$ is an elementary substructure of $\widehat{G_i}$ is a routine consequence of the elimination of quantifiers up to p.p. formulas.  (Note that since $G_i$ is torsion-free and divisible by all primes except for $p_i$, all p.p. formulas are boolean combinations of statements about divisibility by $p_i^k$ and equations, so the fact that $G_i \prec \widehat{G_i}$ follows from the fact that taking the $p_i$-adic completion does not affect divisibility by $p_i$.)  To prove the full theorem, recall that if $K_0 \leq K_1 \leq \ldots$ is an infinite chain of pure extensions then $K_0$ is pure in $\bigcup K_i$.  By this fact plus induction it follows that $H \oplus \bigoplus_{i \in I} G_i$ is a pure subgroup of $H \oplus \bigoplus_{i \in I} \widehat{G_i}$.  To show that the two groups are elementarily equivalent, use the fact that $G_i \equiv \widehat{G_i}$ and Theorem~\ref{direct-sums}.
\end{proof}

\begin{fact}
\label{stability}
(Theorem~3.1 of \cite{prest}) If $R$ is any ring, then any theory in the language of $R$-modules is stable.  In particular, any theory of abelian groups in the pure group language is stable.
\end{fact}

\section{The SB~property and abelian groups}

Our goal in this section is to show that for any abelian group $G$, $\Th(G, +)$ has the SB~property if and only if it is $\omega$-stable (Theorem~\ref{abeliangroups} below).  We warm up by proving a couple of simple lemmas about the SB~property that may be of independent interest, especially Theorem~\ref{pureinj}.

\begin{thm}
\label{inj}
The SB~property holds in the category of all divisible abelian groups, where the arrows consist of \emph{all} group homomorphisms.
\end{thm}

\begin{proof}
Our argument uses the characterization of divisible groups given by Fact~\ref{div_class}.  For any abelian group $G$ and prime $p$, let $G[p] = \left\{g \in G : pg = 0 \right\}$.  Note that if $f: G_0 \rightarrow G_1$ is an injective map between two divisible groups $G_0$ and $G_1$ then for every prime $p$, $f$ maps $G_0 [p]$ into $G_1 [p]$.  The groups $G_0 [p]$ and $G_1 [p]$ are naturally vector spaces over the finite field with $p$ elements and the previous sentence implies that $\dim(G_0 [p]) \leq \dim(G_1 [p])$.  But $\dim(G_i [p])$ is equal to the number of summands of $G_i$ that are isomorphic to $\Z_{p^{\infty}}$, so $G_1$ contains at least as many $\Z_{p^{\infty}}$-summands as $G_0$.  Similarly, if $T_i$ is the torsion subgroup of $G_i$, $f$ induces an injective map from $G_0 / T_0$ into $G_1 / T_1$, so $G_1$ contains at least as many $\Q$-summands as $G_0$.  So if there is also an injective map $g: G_1 \rightarrow G_0$, it follows that $G_0$ and $G_1$ contain the same number of summands of each isomorphism type and hence $G_0 \cong G_1$.
\end{proof}

\begin{thm}
\label{cyc}
The SB~property holds in the category whose objects are all direct sums of finite cyclic groups and whose arrows are all \emph{pure} group embeddings.
\end{thm}

\begin{proof}
The proof is a standard application of Ulm invariants.  To define these, we need some preliminary definitions.  If $p$ is a prime and $G$ is an abelian group, an element $g \in G$ has \emph{$p$-height at least $n$}, or $\textup{ht}_{p}(g) \geq n$, if there is an element $h \in G$ such that $p^n h = g$.  If $G$ is any $p$-group, we can define $P(G, i)$ to be the group of all elements of $G[p]$ of height at least $i$.  The dimension of ${P(G, i)}/{P(G, i + 1)}$ as a vector space over the field with $p$ elements is called the \emph{$i$th Ulm invariant of G}, and we will denote it by $Ulm(G, i)$.  A routine calculation shows that if the $p$-group $G$ is a direct sum of cyclic groups then for any nonzero $i \in \N$, the $i$th Ulm invariant is precisely the number of summands of the form ${\Z}/{p^i \Z}$.

Now suppose that $G_0$ and $G_1$ are two direct sums of finite cyclic groups and $f_0 : G_0 \rightarrow G_1$ and $f_1: G_1 \rightarrow G_0$ are pure injective maps.  It is a standard fact that we can decompose each $G_i$ into a direct sum of $p$-groups (for various primes $p$), so without loss of generality $G_i$ is a $p$-group.  It is routine to check that, for each $k \in \N$ and $i \in \left\{0,1\right\}$, $f_i$ induces a map of vector spaces from ${P(G_i, k)}/{P(G_i, k+1)}$ into ${P(G_{1-i}, k})/{P(G_{1-i}, k+1)}$, and the purity of $f_i$ implies that this induced map is also injective.  Thus the $G_0$ and $G_1$ have the same Ulm invariants.  But this means that for every $k \in \omega$, $G_0$ and $G_1$ contain the same number of summands isomorphic to $\Z / {p^k \Z}$ in any of their direct sum decompositions, so $G_0 \cong G_1$.

\end{proof}

The next result is stated in a little more generality than is needed for the proof of the main theorem.

\begin{thm}
\label{pureinj}
The SB~property holds in the category consisting of all direct sums of indecomposable pure-injective abelian groups (as objects) and all the pure group maps between them (as arrows).
\end{thm}

\begin{proof}
Suppose that $f_0: G_0 \hookrightarrow G_1$ and $f_1: G_1 \hookrightarrow G_0$ be pure embeddings between two such groups.  For $i = 0,1$, write $G_i = K_i \oplus C_i \oplus D_i$ where $D_i$ contains all the divisible summands of $G_i$ (that is, all the summands isomorphic to $\Q$ or to $\Z_{p^{\infty}}$ for some prime $p$), $C_i$ contains all the finite cyclic summands of $G_i$, and $K_i$ contains all the summands isomorphic to $\widehat{\Z_{(p)}}$ for some prime $p$.

\begin{claim}
\label{D-part}
If $x \in G_i$, then $x \in D_i$ if and only if $f_i(x) \in D_{1-i}$.
\end{claim}

\begin{proof}
Any group homomorphism on $G_i$ must map $D_i$ into $D_{1-i}$.  If $f_i(x) \in D_{1-i}$, then the purity of $f_i$ implies that $x \in D_i$.
\end{proof}

\begin{claim}
\label{TD-part}
If $x \in G_i$, then $x \in C_i \oplus D_i$ if and only if $f_i(x) \in C_{1-i} \oplus D_{1-i}$.
\end{claim}

\begin{proof}
If $x = t + d$ where $Nt = 0$ and $d \in D_i$, then $N f_i(t) = f_i(Nt) = 0$ and $f_i(d) \in D_{1-i}$ by purity of $f_i$, so $f(x) = f(t) + f(d) \in C_{1-i} \oplus D_{1-i}$.  Conversely, if $f_i(x) = t' + d'$, $N t' = 0$, and $d' \in D_{1-i}$, then $f(Nx) = N f(x) = N d' \in D_{1-i}$ and thus $Nx \in D_i$.  So we can pick $d \in D_i$ such that $N d = Nx$.  If $T_i$ is the torsion subgroup of $G_i$, then since $N (x - d) = 0$, $x = (x - d) + d \in T_i + D_i$.  But $T_i + D_i = C_i \oplus D_i$, so we are done.
\end{proof}

By Claim~\ref{D-part}, $D_0$ and $D_1$ are purely bi-embeddable, so by Theorem~\ref{inj}, $D_0 \cong D_1$.  By Claim~\ref{D-part} and the ``only if'' part of Claim~\ref{TD-part}, the functions $f_i$ induce well-defined injective maps $\overline{f_i} : {(C_i \oplus D_i)} / {D_i} \hookrightarrow {(C_{1-i} \oplus D_{1-i})} / {D_{1-i}}$.  If $x \in C_i$, $n \in \N$, and $\overline{f_i}(x + D_i) = n (y + D_{1-i}) = ny + D_{1-i}$, then $f_i(x) = ny + d$ for some $d \in D_{1-i}$, so the purity of $f_i$ implies that $x$ is $n$-divisible.  So $x + D_i$ is $n$-divisible in the group ${(C_i \oplus D_i)} / {D_i}$ and $\overline{f_i}$ is a pure embedding.  But ${(C_i \oplus D_i)} / {D_i} \cong C_i$, so $C_0$ and $C_1$ are purely bi-embeddable.  From Theorem~\ref{cyc} we conclude that $C_0 \cong C_1$.

By Claim~\ref{TD-part}, the $f_i$'s also induce well-defined injective maps $\widehat{f_0} : G_0 / {(C_0 \oplus D_0)} \hookrightarrow G_1 / {(C_1 \oplus D_1)}$ and $\widehat{f_1} : G_1 / {(C_1 \oplus D_1)}\hookrightarrow G_0 / {(C_0 \oplus D_0)}$. 

\begin{claim}
The maps $\widehat{f_0}$ and $\widehat{f_1}$ are pure.
\end{claim}

\begin{proof}
We just prove the purity of $\widehat{f_0}$, and the purity of $\widehat{f}_1$ is similar.  If $x \in K_0$, $y \in K_1$, and $\widehat{f}_0 (x + C_0 \oplus D_0) = ny + C_1 \oplus D_1$, then $f_0(x) = ny + t + d$ for some $t \in C_1$ and $d \in D_1$; and if $m \in \N$ is such that $m t = 0$, then $f_0(mx) = mny + md$, so the purity of $f_0$ implies that $mx$ is divisible by $mn$ in $G_0$.  If we pick $z \in G_0$ such that $mx = mnz$, then $m (x - nz) = 0$, so since the torsion subgroup of $G_0$ is contained in $C_0 \oplus D_0$, from $x = nz + (x - nz)$ it follows that $x + C_0 \oplus D_0$ is divisible by $n$ in $G_0 / {(C_0 \oplus D_0)}$.
\end{proof}

From the last claim, it follows that $K_0$ and $K_1$ are purely bi-embeddable.  The $K_i$'s are just direct sums of copies of $\widehat{{\Z}_{(p)}}$ for various primes $p$, and they are determined by the cardinal invariants $\dim(K_i / {p K_i})$ that measure the dimension of the quotient $K_i / {p K_i}$ as a vector space over the finite field with $p$ elements.  Because the $K_i$'s are purely bi-embeddable, $\dim(K_0 / {p K_0}) = \dim(K_1 / {p K_1})$ for every $p$, so $K_0 \cong K_1$.  Since we have also shown that $D_0 \cong D_1$ and $C_0 \cong C_1$, $G_0$ and $G_1$ are isomorphic.

\end{proof}

\begin{definition}
1. $G^\circ$ is the intersection of all $0$-definable subgroups of $G$ of finite index.

2. $\Aut(G / G^\circ)$ is the set of all bijections on $G / G^\circ$ induced by automorphisms of $G$.

3. $\varphi \in \Aut(G / G^\circ)$ is \emph{unipotent} if there is a nonzero $n \in \omega$ such that $\varphi^n$ is the identity map.
\end{definition}

Now we state the main theorem:

\begin{thm}
\label{abeliangroups}
If $G$ is an abelian group, then the following are equivalent:

1. $\Th(G, +)$ has the SB~property;

2. $\Th(G, +)$ is $\omega$-stable;

3. $G$ is the direct sum of a divisible group and a torsion group of bounded exponent;

4. $\Th(G, +)$ is superstable, and if $(\overline{G}, +) \equiv (G, +)$ is saturated, every map in $\Aut(\overline{G} / \overline{G}^\circ)$ is unipotent.
\end{thm}

The proof of Theorem~\ref{abeliangroups} will be proved over the course of several other lemmas and propositions below.  From now on, ``$T$'' will always denote a complete theory of abelian groups in the pure group language.

\begin{lem}
\label{p-adic_NSB}
If $p$ is a prime and $k$ is a positive integer, then $T = \Th(\Z_{(p)}^{k} , +)$ does not have the SB~property.  Furthermore, we can pick a pair of reduced models of $T$ that witness the failure of the SB~property.
\end{lem}

\begin{proof}
For brevity, let $G$ be the additive group of $(\widehat{\Z_{(p)}})^{k}$.  Although we usually consider $\Z$ and $\widehat{\Z_{(p)}}$ as groups, for the purposes of our discussion we will make use of the canonical ring embedding of $\Z$ into $\widehat{\Z_{(p)}}$ so that we can talk about the element ``$1$'' of $\widehat{\Z_{(p)}}$.  For each $i$ between $1$ and $k$, let $e_i$ be the element of $G$ such that the $i$th coordinate of $e_i$ is $1$ and $e_i$'s other coordinates are all $0$.  For any nonempty subset $A \subseteq G$, let $E_G(A)$ ($E$ is for ``envelope'') be the smallest pure subgroup of $G$ containing $A$, that is, if $g \in G, n \in \Z$ and $n g \in E_G(A)$ then $g \in E_G(A)$.  Note that $E_G(A)$ is naturally a $\Z_{(p)}$-module.  If $A$ contains all the $e_i$'s, then since $\Z_{(p)}^k \subseteq E_G(A) \subseteq (\widehat{\Z_{(p)}})^k$ and $\Z_{(p)}^k \prec (\widehat{\Z_{(p)}})^k$ (by Theorem~\ref{completions}), $E_G(A)$ is a model of $T$.

Pick two $p$-adic integers $\gamma_1 , \gamma_2 \in G$ which are units (in the $p$-adic ring) and algebraically independent over $\Z_{(p)}$, that is, if we have $a_{ij} \in \Z_{(p)}$ for $i, j < n$ such that $\Sigma_{i, j < n} a_{ij} {\gamma_1}^{i} {\gamma_2}^{j} = 0$, then each $a_{ij}$ must be zero.  (Such elements exist simply because $\widehat{\Z_{(p)}}$ is an integral domain with uncountably many units, and so by cardinality considerations there must be two units that are algebraically independent over the fraction field of the countable subring $\Z_{(p)}$.)  Let $K_1$ be the subgroup of $G$ generated by $\left\{ {\gamma_1}^i {\gamma_2}^j e_s : 0 \leq i , j < \omega, 1 \leq s \leq k \right\}$.  We like to picture this set as follows:

\begin{equation*}\begin{CD}
e_s @>\gamma_1>> \gamma_1 e_s @>\gamma_1>> \gamma^2_1 e_s @>\gamma_1>> \ldots\\
@VV\gamma_2V @VV\gamma_2V @VV\gamma_2V @.\\
\gamma_2 e_s @>\gamma_1>> \gamma_1 \gamma_2 e_s @>\gamma_1>> \gamma^2_1 \gamma_2 e_s @>\gamma_1>> \ldots\\
@VV\gamma_2V @VV\gamma_2V @VV\gamma_2V @.\\
\gamma^2_2 e_s @>\gamma_1>> \gamma_1 \gamma^2_2 e_s @>\gamma_1>> \gamma^2_1 \gamma^2_2 e_s @>\gamma_1>> \ldots\\
@VV\gamma_2V @VV\gamma_2V @VV\gamma_2V @.\\
\vdots @. \vdots @. \vdots\\
\end{CD}\end{equation*}

Let $K_2$ be the subgroup of $G$ generated by $$\left\{ {\gamma_1}^i {\gamma_2}^j e_s : 1 \leq i < \omega , 0 \leq j < \omega, 1 \leq s \leq k \right\} \cup \left\{e_s : 1 \leq s \leq k \right\},$$ which can be pictured as follows:

\vspace{.2in}

\begin{equation*}\begin{CD}
e_s @>\gamma_1>> \gamma_1 e_s @>\gamma_1>> \gamma^2_1 e_s @>\gamma_1>> \ldots\\
@. @VV\gamma_2V @VV\gamma_2V @.\\
 @. \gamma_1 \gamma_2 e_s @>\gamma_1>> \gamma^2_1 \gamma_2 e_s @>\gamma_1>> \ldots\\
@. @VV\gamma_2V @VV\gamma_2V @.\\
@. \gamma_1 \gamma^2_2 e_s @>\gamma_1>> \gamma^2_1 \gamma^2_2 e_s @>\gamma_1>> \ldots\\
@. @VV\gamma_2V @VV\gamma_2V @.\\
 @. \vdots @. \vdots\\
\end{CD}\end{equation*}

Let $H_1 := E_G(K_1)$ and $H_2 := E_G(K_2)$.  Then both $H_1$ and $H_2$ are elementary submodels of $G$ by the comment at the end of the previous paragraph and $H_1$ and $H_2$ are both reduced groups.

Note that $G$ is not only a $\Z_{(p)}$-module but also a $\widehat{\Z_{(p)}}$-module (using ring multiplication in $\widehat{\Z_{(p)}}$ to define the action $\alpha \cdot (a_1, \ldots, a_k) = (\alpha a_1, \ldots, \alpha a_k)$.)  If $\alpha \in \widehat{\Z_{(p)}}$ is a ring unit, let $\sigma_{\alpha}$ be the group automorphism of $G$ given by $\sigma_{\alpha}(x) = \alpha \cdot x$.  

\begin{claim}
$\sigma_{\gamma_1}(H_1) \subseteq H_2$.
\end{claim}

\begin{proof}
First notice that $\sigma_{\gamma_1} (K_1) \subseteq K_2$.  As a group, $H_1$ is generated by elements of the form $n^{-1} a$, where $n \in \N$ and $a \in K_1$.  (We are using the fact that since $G$ is torsion-free, if $a \in G$ and there is an element $b \in G$ such that $n b = a$, then there is only one such $b$, so we can write ``$b = n^{-1} a$.'')  So it suffices to show that if $a \in K_1$ is $n$-divisible in $G$ then $\sigma_{\gamma_1}(n^{-1} a) \in H_2$.   But since $n \cdot \sigma_{\gamma_1}(n^{-1} a) = \sigma_{\gamma_1}(n n^{-1} a) = \sigma_{\gamma_1}(a)$, $\sigma_{\gamma_1}(n^{-1} a) = n^{-1} \sigma_{\gamma_1}(a) \in E_G(K_2) = H_2$.
\end{proof}

Since $H_2 \prec G$, $\sigma_{\gamma_1}(H_1) \prec G$, and $\sigma_{\gamma_1}(H_1) \subseteq H_2$, it follows that $\sigma_{\gamma_1}(H_1) \prec H_2$.  Since $H_2 \subseteq H_1$, the same argument shows that $H_2 \prec H_1$, so the $H_i$'s are elementarily bi-embeddable and the only thing left to check is that $H_1$ and $H_2$ are not isomorphic as groups.

\begin{claim}
\label{mult}
If $H \leq G$ is a subgroup containing $e_1, \ldots, e_k$ and $f: H \rightarrow G$ is any elementary map, then there is a unique invertible $k \times k$ matrix $A$ over $\widehat{\Z_{(p)}}$ such that for all $x \in H$, $f(x) = A x$.  Consequently, $f$ has a unique extension to an automorphism of $G$.
\end{claim}

\begin{proof}
By elementarity, for each $n \in \N$, $f$ induces a well-defined group automorphism $f_n$ of ${G}/{p^n G}$.  Each automorphism $f_n$ is additive, so $f_n$ is determined uniquely by the $k$ values $f_n(e_1 + p^n G), \ldots, f_n(e_k + p^n G)$ and $f_n$ can be represented as an invertible $k \times k$ matrix $A_n$ over ${\widehat{\Z_{(p)}}}/{p^n \widehat{\Z_{(p)}}}$ whose columns list the $k$ coordinates of the elements $f_n(e_i + p^n G)$.  The matrix $A = \lim_{n < \omega} A_n$ is what we are looking for.  If $B_n$ is the $k \times k$ matrix over ${\widehat{\Z_{(p)}}}/{p^n \widehat{\Z_{(p)}}}$ such that $A_n  B_n = I$ and $B = \lim_{n < \omega} B_n$, then $A B = \lim_{n < \omega} A_n \cdot \lim_{n < \omega} B_n = \lim_{n < \omega} A_n B_n = \lim_{n < \omega} I = I$.
\end{proof}

For $\ell = 1, 2$, let $S_\ell$ be the subring of the \emph{group} endomorphism ring of $H_\ell$ which is generated by all pure and injective endomorphisms of $H_\ell$ and let $R_\ell$ be the center of $S_\ell$.  By Claim~\ref{mult}, the $R_\ell$'s are isomorphic to subrings of the $k \times k$ matrix ring $M_{kk} (\widehat{\Z_{(p)}})$.

\begin{claim}
\label{centers}
For $\ell = 1,2$, every element of $R_\ell$ is represented by a diagonal matrix all of whose diagonal entries are equal.
\end{claim}

\begin{proof}
Let $e_{ij}$ be the $k \times k$ matrix whose $(i,j)$ entry is $1$ and all of whose other entries are $0$.  Since the matrix $I + e_{ij}$ induces an automorphism of each $H_{\ell}$, $e_{ij} = (I + e_{ij}) - I$ is in $S_\ell$.  If $A \in S_\ell$ has a nonzero off-diagonal entry $a_{ij}$, then the usual computation shows that $A e_{ji} \neq e_{ji} A$, so every element of $R_\ell$ is a diagonal matrix.  Since every permutation matrix is in $R_\ell$ (since permutation of coordinates is an automorphsim of $H_{\ell}$), all the diagonal entries of any matrix in $R_\ell$ must be equal.
\end{proof}

By the last claim, each map $r \in R_\ell$ is equal to $\sigma_{\alpha}$ for some $\alpha \in \widehat{\Z_{(p)}}$.  We show that $R_1$ is not isomorphic to $R_2$, which is enough to show that $G_1$ is not isomorphic to $G_2$ since the $R_\ell$'s are defined without reference to the embeddings of the $H_\ell$'s into $G$.

More notation: if $R$ is any integral domain of characteristic zero (that is, the field of fractions of $R$ has characteristic zero) and $A \subset R$, then let $E_R(A)$ be the smallest subring of $R$ containing $A$ with the property that if $n \in \Z$, $r \in R$, and $n r \in E_R(A)$, then $r \in E_R(A)$.  Notice that any element of $E_R(A)$ is a $\Q$-linear combination of elements of $A$  Therefore for any $r \in E_R(A)$ there is a finite subset $A_0 \subseteq A$ such that $r \in E_R(A_0)$.

\begin{claim}
$R_1 = E_{R_1}(\sigma_{\gamma_1} , \sigma_{\gamma_2})$ and $R_2 = E_{R_2}(\sigma_{\gamma_1 \gamma_2^i} : 0 \leq i < \omega)$.
\end{claim}

\begin{proof}
It is routine to check that $\sigma_{\gamma_1}, \sigma_{\gamma_2} \in R_1$.  Suppose that $r \in R_1$ and $r = \sigma_{\alpha}$ for $\alpha \in \widehat{\Z_{(p)}}$.  Then $\sigma_{\alpha}(e_1) \in R_1$ implies that $\alpha$ is some $\Q$-linear combination of the elements $\left\{\gamma_1^i \gamma_2^j : i, j \in \N\right\}$, and hence $\sigma_{\alpha}$ is a $\Q$-linear combination of the elements $\left\{\sigma_{\gamma_1^i \gamma_2^j} : i, j \in \N\right\}$.  But $\sigma_{\gamma_1^i \gamma_2^j} = (\sigma_{\gamma_1})^i(\sigma_{\gamma_2})^j$ is in $E_{R_1}(\sigma_{\gamma_1}, \sigma_{\gamma_2})$ and $E_{R_1}(\sigma_{\gamma_1}, \sigma_{\gamma_2})$ is closed under $\Q$-linear combinations inside $R_1$, so $\sigma_{\alpha} \in E_{R_1}(\sigma_{\gamma_1}, \sigma_{\gamma_2})$.

The case of $R_2$ is similar.
\end{proof}

The $\Q$-independence of $\gamma_1$ and $\gamma_2$ implies that $E_{R_2}(\sigma_{\gamma_1}) \subset E_{R_2}(\sigma_{\gamma_1}, \sigma_{\gamma_1 \gamma_2}) \subset E_{R_2}(\sigma_{\gamma_1}, \sigma_{\gamma_1 \gamma_2}, \sigma_{\gamma_1 \gamma_2^2}) \subset \ldots$ is a strictly-ascending tower of subrings of $R_2$ whose union is the whole ring.  Thus there cannot exist two elements $\alpha , \beta \in R_2$ such that $R_2 = E_{R_2}(\alpha, \beta)$ and $R_1 \ncong R_2$.

\end{proof}

\begin{cor}
\label{p-adic_NSB2}
Suppose that $\langle p_i : i \in I \rangle$ is a sequence of distinct primes and $\langle k_i : i \in I \rangle$ is a sequence of positive integers.  Then $T = \Th(\bigoplus_{i \in I} \Z_{(p_i)}^{k_i}, +)$ does not have the SB~property.
\end{cor}

\begin{proof}
Suppose that for each $i \in I$ we have two reduced models $H_{i, 1}$ and $H_{i,2}$ of $\Th(\Z_{(p_i)}^{k_i})$ and $\varphi$ is a group map from $\bigoplus_{ i \in I} H_{i,1}$ into $\bigoplus_{i \in I} H_{i,2}$.  Then if $h \in H_{i,1}$, $\varphi(h)$ must be in $H_{i,2}$, since if $i \neq j$ and $\varphi(h)$ has a nonzero projection $a$ onto the component $H_{j, 2}$, $a$ is a nonzero element of $H_{j,2}$ with infinite $p_{j}$-height, contradicting the fact that $H_{j,2}$ is reduced.  Since the counterexamples to the SB~property constructed in Lemma~\ref{p-adic_NSB} are reduced, it follows that we can get a counterexample to the SB~property for $T$ by taking direct sums of counterexamples of the SB~property for each of the theories $\Th(\Z_{(p_j)}^{k_j})$.

\end{proof}

\begin{cor}
\label{p-adic_NSB3}
If $T$ has a pure-injective model of the form $H = K \oplus C \oplus D$, where $K$ is a nonzero direct sum of copies of $\widehat{\Z_{(p)}}$ for various primes $p$, $C$ is a direct sum of finite cyclic groups, and $D$ is the maximal divisible subgroup, then $T$ does not have the SB~property.
\end{cor}

\begin{proof}
By Corollary~\ref{p-adic_NSB2}, there are models $K_0, K_1$ of $\Th(K, +)$ that are purely bi-embeddable but not isomorphic.  By Theorem~\ref{completions}, $G_0 := K_0 \oplus C \oplus D$ and $G_1 := K_1 \oplus C \oplus D$ are both models of $T$, and $G_0$ and $G_1$ are purely bi-embeddable.  Any group isomorphism $\varphi : G_0 \rightarrow G_1$ would have to map $D$ onto itself, and so it would induce a group isomorphism $\varphi' : K_0 \oplus C \rightarrow K_1 \oplus C$; and since $K_0$ and $K_1$ are torsion-free, $\varphi'$ would map $H$ onto itself, so there would be a group isomorphism from $K_0$ onto $K_1$, contradiction.  Thus $G_0 \ncong G_1$.
\end{proof}

\begin{thm}
\label{nice_case}
Let $S$ be an infinite set of prime numbers and let $\langle r_p : p \in S \rangle$ be a sequence of positive integers.  Then the group $\bigoplus_{p \in S} \left(\Z / {p \Z} \right)^{r_p}$ does not have the SB~property.  Moreover, the fact that this theory does not have the SB~property can be witnessed by a pair of reduced groups.
\end{thm}

\begin{proof}
Some notation: let $G = \bigoplus_{p \in S} \left(\Z / {p \Z} \right)^{r_p}$ and let $\widehat{G} =  \prod_{p \in S} \left(\Z / {p \Z} \right)^{r_p}$.  Notice that $\End(G) \cong \prod_{p \in S} \End(\left(\Z / {p \Z} \right)^{r_p})$, and let $R$ be the subring of $\End(G)$ consisting of all endomorphisms that correspond to elements of the form $\prod_{p \in S} \alpha_p \in \prod_{p \in S} \End(\left(\Z / {p \Z} \right)^{r_p})$ where each $\alpha_p \in \End(\left(\Z / {p \Z} \right)^{r_p})$ is multiplication by some element of $\Z$.  Then $R$ is commutative.  Let $A \subseteq R$ be the subset consisting of all elements of $R$ that are automorphisms of $G$.  ($A$ is not a subring since it does not contain $0$ and is not closed under addition.)  Observe that $R$ is a Polish space under the product topology and that $A$, being a closed subset of $A$, is also a Polish space.  For each prime $p \in \N$, we define a map ``$\alpha_{\frac{1}{p}}$'' on $\widehat{G}$ as follows: if $g \in \widehat{G}$ is divisible by $p$, then $\alpha_{\frac{1}{p}}( g)$ is the unique element $h \in \widehat{G}$ such that $p h = g$ and $h$ is also divisible by $p$; and if $g \in \widehat{G}$ is not divisible by $p$, let $h \in \left(\Z / {p \Z} \right)^{r_p}$ be the unique element such that $g + h$ is divisible by $p$, and we define $\alpha_{\frac{1}{p}}(g) = \alpha_{\frac{1}{p}}(g + h)$.  Note that this map $\alpha_{\frac{1}{p}}$ is in $R$.  If $p$ is a prime that is not in $S$, then we define $\alpha_{\frac{1}{p}}$ to be the map that takes every $g \in \widehat{G}$ to the unique $h \in \widehat{G}$ such that $p h = g$.  Similarly, for any $n \in \N$, we can define a map $\alpha_{\frac{1}{n}}$ using the definitions of the $\alpha_{\frac{1}{p}}$'s and the rule $\alpha_{\frac{1}{ab}} = \alpha_{\frac{1}{a}} \circ \alpha_{\frac{1}{b}}$.

Let $F \subset R$ be the subset of $R$ consisting of all elements of the form $\prod_{p \in S} \alpha_p$ where almost every $\alpha_p$ is the zero map.

\begin{claim}
For any nonzero $q \in \Z[x,y]$, $n \in \N$, and $\alpha \in F$, the set $$X_{q, n, \alpha} = \left\{(\sigma,\tau) \in A \times A : \alpha_{\frac{1}{n}} \left[ q(\sigma,\tau)\right] = \alpha \right\}$$ is closed and nowhere dense.
\end{claim}

\begin{proof}
First, note that in order for $X_{q, n, \alpha}$ to be nonempty, it is necessary that for every prime $p \in S$ that divides $n$, the projection of $\alpha$ onto $\left(\Z / {p \Z}\right)^{r_p}$ is zero.  So from now on we will assume that this is true and consider only the set $S' \subseteq S$ consisting of those primes in $S$ that do not divide $n$.  If $(\sigma,\tau) \in A^2$, $(\sigma,\tau) \in X_{q,n,\alpha}$ if and only if for every $p \in S'$, the projections $\sigma_p$ and $\tau_p$ of $\sigma$ and $\tau$ onto $\left(\Z / {p \Z}\right)^{r_p}$ satisfy a certain polynomial $q_p(x,y)$.  Suppose that $\langle (\sigma(i), \tau(i)) : i \in \N \rangle$ is some convergent sequence in $A^2$ such that each term $(\sigma(i), \tau(i))$ is in $X_{q,n,\alpha}$.  Then for every $p \in S'$, $(\sigma(i)_p, \tau(i)_p)$ satisfies $q_p(x,y)$.  So since the convergence of the sequence implies that for each $p \in S'$, $(\sigma(i)_p, \tau(i)_p)$ is eventually constant, its limit $(\sigma_p, \tau_p)$ also satisfies $q_p(x,y)$.  But this means that the limit $(\sigma,\tau)$ of the sequence $\langle (\sigma(i), \tau(i)) : i \in \omega \rangle$ is in the set $X_{q,n,\alpha}$, so $X_{q,n,\alpha}$ is closed.

Since $X_{q,n,\alpha}$ is closed, to show that it is nowhere dense it suffices to show that $X_{q,n,\alpha}$ contains no basic open subsets.  Suppose that $U$ is a basic open subset of $A^2$ defined by $$U = \left\{(\sigma,\tau) \in A^2 : \textup{ for every } p \in T, \, (\sigma_p, \tau_p) \in Y_p \right\},$$ where $T$ is some finite subset of $S$ and $Y_p$ is a subset of $\left(\Z / {p \Z}\right)^{r_p} \times \left(\Z / {p \Z}\right)^{r_p}$.  Then if $S'' = S' - (T \cup \left\{p \in S : \alpha_p \neq 0 \right\})$, $S''$ is an infinite set of primes such that for every $p \in S''$, the evaluation of $q(x,y)$ is the zero function from $\mathbf{F}_p \times \mathbf{F}_p$ into $\mathbf{F}_p$.  But this contradicts the fact that $q$ is nonzero, since if $m$ is the number of nonzero terms in $q$ and $p > m$ then $q(1,1)$ is not zero in $\mathbf{F}_p$.
\end{proof}

By the claim, the countability of $F$, and the Baire category theorem, we can pick $(\sigma_1, \sigma_2) \in A \times A$ such that $(\sigma_1, \sigma_2)$ does not lie in any of the sets $X_{q, n, \alpha}$ for nonconstant $q \in \Z[x,y]$.  Let $a \in \widehat{G}$ be some element whose projection onto every $\Z / {p \Z}$ (for $p \in S$) is nonzero.  Let $K_1$ be the set $\left\{\sigma_1^i \sigma_2^j (a) : i, j \in \omega \right\}$, which can be pictured as follows:

\begin{equation*}\begin{CD}
a @>\sigma_1>> \sigma_1 (a) @>\sigma_1>> \sigma^2_1 (a) @>\sigma_1>> \ldots\\
@VV\sigma_2V @VV\sigma_2V @VV\sigma_2V @.\\
\sigma_2 (a) @>\sigma_1>> \sigma_1 \sigma_2 (a)@>\sigma_1>> \sigma^2_1 \sigma_2 (a) @>\sigma_1>> \ldots\\
@VV\sigma_2V @VV\sigma_2V @VV\sigma_2V @.\\
\sigma^2_2 (a) @>\sigma_1>> \sigma_1 \sigma^2_2 (a)@>\sigma_1>> \sigma^2_1 \sigma^2_2 (a) @>\sigma_1>> \ldots\\
@VV\sigma_2V @VV\sigma_2V @VV\sigma_2V @.\\
\vdots @. \vdots @. \vdots\\
\end{CD}\end{equation*}

Let $K_2$ be the set $\left\{a \right\} \cup \left\{\sigma_1^i \sigma_2^j (a) : 0 < i < \omega, j < \omega \right\}$, as in the diagram below:

\begin{equation*}\begin{CD}
a @>\sigma_1>> \sigma_1 (a) @>\sigma_1>> \sigma^2_1 (a) @>\sigma_1>> \ldots\\
@. @VV\sigma_2V @VV\sigma_2V @.\\
 @. \sigma_1 \sigma_2 (a) @>\sigma_1>> \sigma^2_1 \sigma_2 (a)@>\sigma_1>> \ldots\\
@. @VV\sigma_2V @VV\sigma_2V @.\\
@. \sigma_1 \sigma^2_2 (a)@>\sigma_1>> \sigma^2_1 \sigma^2_2 (a)@>\sigma_1>> \ldots\\
@. @VV\sigma_2V @VV\sigma_2V @.\\
 @. \vdots @. \vdots\\
\end{CD}\end{equation*}

For $\ell = 1,2$, let $H_\ell$ be the smallest pure subgroup of $\widehat{G}$ containing $G \cup K_\ell$.

\begin{claim}
\label{h-ell-rep}
$H_\ell$ is the set of all finite sums of elements of the form $\alpha_{\frac{1}{n}}(g + k)$ where $n \in \N$, $g \in G$, and $k$ is a finite $\Z$-linear combination of elements of $K_\ell$.
\end{claim}

\begin{proof}
Call $Z$ the set of all such finite sums.  Since $\alpha_{\frac{1}{n}}(-g -k) = - \alpha_{\frac{1}{n}}(g + k)$, it is clear that $Z$ is a group.  Suppose that the element $a = \sum_{1 \leq i \leq s} \alpha_{\frac{1}{n_i}}(g_i + k_i)$ of $Z$ is divisible by $m$ in $\widehat{G}$.  Let $S_0 \subseteq S$ be the set of all primes in $S$ that occur in the prime factorization of $m$ and let $G_0 = \sum_{p \in S_0} \left(\Z / {p \Z} \right)^{r_p}$.  If $h_i$ is the projection of $\alpha_{\frac{1}{n_i}}(g_i + k_i)$ onto $G_0$, then $\sum_{1 \leq i \leq s} h_i = 0$, so if $h'_i = \sum_{1 \leq j \leq s, j \neq i} h_j$, we have that each $\alpha_{\frac{1}{n_i}}(g_i + k_i) + h'_i$ is divisible by $m$ and $\sum_{1 \leq i \leq s}h'_i = 0$.

Thus $$\alpha_{\frac{1}{m}}(a) = \alpha_{\frac{1}{m}}\left(\sum_{1 \leq i \leq s} \alpha_{\frac{1}{n_i}}(g_i + k_i) \right) = \alpha_{\frac{1}{m}}\left(\sum_{1 \leq i \leq s} \alpha_{\frac{1}{n_i}}(g_i + k_i) + h'_i \right)$$ $$= \sum_{1 \leq i \leq s} \alpha_{\frac{1}{m}}\left(\alpha_{\frac{1}{n_i}}(g_i + k_i) + h'_i \right) = \sum_{1 \leq i \leq s} \alpha_{\frac{1}{m}}\left(\alpha_{\frac{1}{n_i}}(g_i + k_i) \right) = \sum_{1 \leq i \leq s} \alpha_{\frac{1}{m n_i}}(g_i + k_i),$$

which is also in $Z$.  So $Z$ is a pure subgroup of $\widehat{G}$.  Since $Z \subseteq H_\ell$, minimality implies that $Z = H_\ell$.
\end{proof}

Let $R_\ell$ be the ring $\left\{ r \in R: r(H_\ell) \subseteq H_\ell \right\}$.  If $a, b \in R$, then we write $a \sim b$ if $a -b \in F$.  Note that since $G \subseteq H_\ell$, $F \subseteq R_\ell$, and so if $a \in R_\ell$, $b \in R$, and $a \sim b$, then $b \in R_\ell$.  If $B \subseteq R_{\ell}$, then we define $E_{R_\ell}(B)$ to be the smallest subring of $R_{\ell}$ containing $B \cup F \cup \left\{\alpha_{\frac{1}{n}} : n \in \omega \right\}$.

\begin{claim}
\label{R2-chain}
$R_1 = E_{R_1}(\sigma_1, \sigma_2)$ and $R_2 = E_{R_2}(\sigma_1, \sigma_1 \sigma_2, \sigma_1  \sigma_2^2, \ldots)$.
\end{claim}

\begin{proof}
Since $\sigma_1, \sigma_2 \in R_1$, it is clear that $E_{R_1}(\sigma_1, \sigma_2) \subseteq R_1$.  For the other direction, suppose that $\tau \in R_1$.  Since $\tau(a) \in R_1$, by Claim~\ref{h-ell-rep} this element can be written as $$\tau(a) = \sum_{1 \leq i \leq s} \alpha_{\frac{1}{n_i}} (g_i + k_i),$$ where the $n_i$'s are in $\N$, the $g_i$'s are in $G$, and the $k_i$'s are finite $\Z$-linear combinations of elements in $K_1$.  Since each element $\alpha_{\frac{1}{n_i}}(g_i + k_i) - \alpha_{\frac{1}{n_i}}(k_i)$ is in $G$, it follows that $\tau(a) - \sum_{1 \leq i \leq s} \alpha_{\frac{1}{n_i}} (k_i) \in G$.  Notice that there is a unique $\tau_i \in R_1$ such that $\tau_i (a) = k_i$, and $\tau_i = q_i(\sigma_1, \sigma_2)$ for some polynomial $q_i(x,y) \in \Z[x,y]$.  So on all but finitely many coordinates, $\tau$ is equal to $\sum_{1 \leq i \leq s} \alpha_{\frac{1}{n_i}}\left[q_i(\sigma_1, \sigma_2)\right] \in E_{R_1}(\sigma_1, \sigma_2)$, and so $\tau$ itself is in $E_{R_1}(\sigma_1, \sigma_2)$.  The proof for $R_2$ is almost identical.
\end{proof}

\begin{claim}
\label{prop-incl}
For every $i \in \N$, $\sigma_1 \sigma_2^{m+1} \notin E_{R_2}(\sigma_1, \sigma_1 \sigma_2, \ldots, \sigma_1 \sigma_2^m)$.
\end{claim}

\begin{proof}
Suppose towards a contradiction that $\sigma_1 \sigma_2^{m+1} \in E_{R_2}(\sigma_1, \sigma_1 \sigma_2, \ldots, \sigma_1 \sigma_2^m)$, so that $\sigma_1 \sigma_2^{m+1} \sim \tau$ for some $\tau$ in the subring $R_{2,m} \subseteq R_2$ generated by the set $\left\{\sigma_1, \ldots, \sigma_1 \sigma_2^m \right\} \cup \left\{\alpha_{\frac{1}{n}} : n \in \N \right\}$.  Since $\tau \in R_{2,m}$, $\tau \sim \alpha_{\frac{1}{n}} \left[q(\sigma_1, \sigma_2) \right]$ for some $n \in \N$ and some $q \in \Z[x,y]$ such that for every nonzero term $a x^i y^j$  of $q$, if $j > m$ then $i > 1$.  This condition on $q$ implies that $q - xy^{m+1} \neq 0$.  But $\sigma_1 \sigma_2^{m+1} \sim \alpha_{\frac{1}{n}} \left[q(\sigma_1, \sigma_2)\right]$ implies that there is some $\tau_0 \in F$ such that $\tau_0 \sim \alpha_{\frac{1}{n}} \left[q(\sigma_1, \sigma_2) \right] - \sigma_1 \sigma_2^{m+1}  \sim \alpha_{\frac{1}{n}} \left[(q - x y^{m+1})(\sigma_1, \sigma_2) \right]$, which contradicts how $\sigma_1$ and $\sigma_2$ were chosen.
\end{proof}

Now if $H_1$ and $H_2$ were isomorphic, then there would be an isomorphism $\varphi: R_1 \rightarrow R_2$, and necessarily $\varphi$ would map the set $F$ onto itself and fix $\left\{\alpha_{\frac{1}{n}} : n \in \omega \right\}$ pointwise.  This means that if $A \cup \left\{ \tau \right\} \subseteq R_1$, then $\tau \in E_{R_1} (A)$ if and only if $\varphi(\tau) \in E_{R_2}(\varphi(A))$.  This yields a contradiction, since $R_1 = E_{R_1}(\sigma_1, \sigma_2)$ but Claims~\ref{R2-chain} and \ref{prop-incl} imply that there is no finite subset $B \subseteq R_2$ such that $R_2 = E_{R_2}(B)$.

For the ``Moreover...'' clause, note that our groups $H_1$ and $H_2$ are reduced since they are subgroups of the reduced group $\widehat{G}$.

\end{proof}

\begin{thm}
\label{unbdd_tor_SB}
If $G$ is a direct sum of cyclic groups of unbounded order -- that is, there is no $n \in \N$ such that $nG = 0$ -- then $\Th(G, +)$ does not have the SB~property.  Moreover, the fact that $\Th(G, +)$ does not have the SB~property can be witnessed by a pair of reduced groups.
\end{thm}

\begin{proof}

Suppose $G = \bigoplus_{k \in \omega, p \textup{ prime}} (\Z / {p^k \Z})^{\kappa_{k, p}}$.

\begin{claim}
\label{primebound}
We may assume that for each prime $p$, there is a number $m_p \in \N$ such that for every $k > m_p$, $\kappa_{k, p} = 0$.
\end{claim}

\begin{proof}
If there is a prime $p$ for which this is not true, then for every $k \in N$, $p^{k+1} G$ has infinite index in $p^k G$, so $G$ is not even superstable, and Theorem~5.5 of \cite{when} implies that $G$ does not have the SB~property.
\end{proof}

\begin{claim}
We may assume that every $\kappa_{k, p}$ is finite.
\end{claim}

\begin{proof}
If infinitely many of the $\kappa_{k,p}$'s are infinite, then as in the previous claim, $G$ is not even supertable and thus $G$ does not have the SB~property.  If only finitely many of the $\kappa_{k,p}$'s are infinite, let $(k_1, p_1), \ldots, (k_\ell, p_\ell)$ be a list of all the pairs such that $\kappa_{k_i, p_i}$ is infinite.  By Claim~\ref{primebound}, for each $i \leq \ell$ there is some $k'_i \in \N$ such that for every $j > k'_i$, $\kappa_{j, p_i} = 0$.  Let $M = p_1^{k'_1} \cdot \ldots \cdot p_{\ell}^{k'_{\ell}}$.  Then we can write $G = G[M] \oplus M G$, where $G[M] = \left\{g \in G : Mg = 0 \right\}$.
So if $\Th(M G)$ does not have the SB~property, as witnessed by the two bi-embeddable, nonisomorphic groups $K_1$ and $K_2$, then Fact~\ref{direct-sums} implies that $G [M] \oplus K_1$ and $G [M] \oplus K_2$ are purely bi-embeddable models of $\Th(G, +)$, and they are nonisomorphic because any isomorphism between them would induce an isomorphism between $K_1$ and $K_2$.

\end{proof}

\begin{claim}
We may assume that for every $p$ and every $k > 1$, $\kappa_{k,p} = 0$.
\end{claim}

\begin{proof}

Recall that the \emph{socle of $G$} is the direct sum of all the minimal subgroups of $G$ (where $H \leq G$ is \emph{minimal} if $H \neq 0$ and $H$ has no nontrivial subgroups -- so this is not the model-theoretic sense of ``minimal.'').  If $H$ is the socle of $G$, then there is some infinite set $S$ of prime numbers and some sequence $\langle r_p : p \in S \rangle$ of numbers $r_p \in \N$ such that $H \cong \bigoplus_{p \in S} \left(\Z / {p \Z} \right)^{r_p}$.

Suppose that $K \cong H$.  Then if $k \in K$, $k \neq 0$, $p \in S$, and $p$ divides $k$, then $p^n k \neq 0$ for every $n \in \N$.  This implies that the maximal divisible subgroup of $K$ is torsion-free, so by Fact~\ref{div_class}, $K \cong R \oplus {\Q}^{\lambda}$ where $R$ is some reduced subgroup of $K$.  Note that the torsion subgroup of $R$ is isomorphic to $H$.  Also, if $r_1, r_2 \in R$ and there is a sequence $\langle a_p : p \in S \rangle$ of elements in $R$ such that $p a_p = 0$ and $r_\ell - a_p$ is divisible by $p$ (for $\ell = 1, 2$), then $r_1 = r_2$, from which it follows that $R$ is isomorphic to a subgroup of $\prod_{p \in S} \left(\Z / {p \Z} \right)^{r_p}$.

Now suppose that $H$ does not have the SB~property, witnessed by two bi-embeddable, nonisomorphic models $H_1$ and $H_2$.  Then by the previous paragraph the divisible parts of $H_1$ and $H_2$ are isomorphic and we can assume without loss of generality that both $H_1$ and $H_2$ are subgroups of $\prod_{p \in S} \left(\Z / {p \Z} \right)^{r_p}$.  For $\ell = 1, 2$, there is a unique pure subgroup $K_\ell$ of $\prod_{k \in \omega, p \textup{ prime}} (\Z / {p^k \Z})^{\kappa_{k, p}}$ such that $G$ is a subgroup of $K_\ell$, the socle of $K_\ell$ is $H_\ell$, and $K_\ell$ is a subgroup of any other group with all these properties.  Then $K_1$ and $K_2$ are purely bi-embeddable and nonisomorphic (since their socles are nonisomorphic), so $\Th(G, +)$ does not have the SB~propery.

\end{proof}

To complete the proof of Theorem~\ref{unbdd_tor_SB}, apply the last three Claims and Theorem~\ref{nice_case}.

\end{proof}

\emph{Proof of Theorem~\ref{abeliangroups}:}

By Theorem~\ref{cases}, $\Th(G, +)$ has a model that is a direct sum of indecomposable pure-injective groups.  Let $H$ be such a model.

1 $\Rightarrow$ 3:  Suppose that $T$ has the SB~property.

Case~A: For some prime $p$, there are infinitely many $k \in \omega$ such that $\Z / {p^k \Z}$ is a direct summand of $H$.

In Case~A, the chain $H \supseteq p H \supseteq p^2 H \supseteq \ldots$ of definable subgroups has the property that for any $k \in \N$, $p^{k+1} H$ has infinite index in $p^k H$, so $T$ is not superstable, contradicting Theorem~5.5 of \cite{when}.  So Case~A cannot hold.

Case~B: For infinitely many primes $p$, there is a $k \in \omega$ such that $\Z / {p^k \Z}$ is a direct summand of $H$.

In Case~B, let $H = A \oplus B$ where $A$ is the torsion subgroup of $H$.  By Theorem~\ref{unbdd_tor_SB}, there are two reduced, bi-embeddable, nonisomorphic models $A_1$ and $A_2$ of $\Th(A, +)$.  For $\ell = 1$ or $2$, let $H_\ell = A_\ell \oplus B$.  Then $H_1$ and $H_2$ are bi-embeddable models of $\Th(H, +)$.  To see that they are isomorphic, let $C_\ell$ be the set of all $g \in H_\ell$ such that $g$ is \emph{not} in the divisible subgroup of $H_\ell$ and there exists a sequence of torsion elements $\langle a_i : i \in \omega \rangle$ such that for every $i \in \omega$, $g - a_i$ is divisible by $n$ in $H_\ell$.  Note that $C_\ell$ is defined independently of how we choose to represent $H_\ell$ as a direct sum of $A_\ell$ and pure-injective groups, and $A_\ell = C_\ell \cup \left\{0\right\}$, so the isomorphism class of $A_\ell$ is an invariant of the isomorphism class of $H_\ell$.

So we are not in Case~A or Case~B and there are primes $p_1, \ldots, p_n$ and positive integers $k_1, \ldots, k_n$ such that $H$ can be written as $$\bigoplus_{i \leq n} \left(\Z / {p^{k_i} \Z}\right)^{\kappa_i} \oplus \bigoplus_{p \textup{ prime}} \left(\Z_{p^{\infty}}\right)^{\lambda_p} \oplus \Q^{\mu}, \, \, \, \, \,(*)$$ where $\kappa_i$, $\lambda_p$, and $\mu$ are cardinal numbers.  Therefore any group $K$ that models $\Th(G, +)$ must have the property that there is a number $m \in \N$ such that $m K$ is divisible.  From Fact~\ref{div-summand}, Fact~\ref{div_class}, and the fact that any abelian group $A$ such that $m A = 0$ is a direct sum of finite cyclic groups (see Section~11 of \cite{kaplansky}), it follows that any model of $\Th(G, +)$ is of the form (*), so in particular $G$ itself must be a direct sum of a divisible group and cyclic groups of bounded order.

3 $\Rightarrow$ 2:
As in the previous paragraph, any model of $\Th(G, +)$ must be of the form (*), and $\omega$-stability of $T$ follows from the elimination of quantifiers up to p.p. formulas and a standard type-counting argument.

2 $\Rightarrow$ 1: We have already shown that in Case~A above, the theory of the group is not $\omega$-stable (since in Case~A the theory is not even superstable).  In Case~B, if we use the L\"{o}wenheim-Skolem theorem to get a countable submodel $H_0 \prec H$, then $H_0 \oplus \Z_{(p)} \prec H \oplus \widehat{\Z_{(p)}}$.  Any two distinct elements of $\widehat{\Z_{(p)}}$ realize different types over the countable model $H_0 \oplus \Z_{(p)}$, so again the theory is not $\omega$-stable.  Thus if $T$ is $\omega$-stable then the argument from 1 $\Rightarrow$ 2 shows that any model of $T$ must be of the form (*).  This means that any model of $T$ is a direct sum of indecomposable pure-injective groups, and the SB~property follows from Theorem~\ref{pureinj} above.

2 $\Rightarrow$ 4: If $G$ is $\omega$-stable, then $\left[G : G^\circ \right]$ is finite.  (This is by standard arguments: otherwise, $\left[G : G^\circ\right] = 2^{\aleph_0}$, and over any countable model of $\Th(G, +)$ we could define uncountably many distinct cosets of $G^\circ$, contradicting $\omega$-stability.)  So clearly any $\varphi \in \Aut(G / G^{\circ})$ is unipotent.

4 $\Rightarrow$ 3: We claim that 4 implies that no model $K$ of $\Th(G, +)$ can be of the form $H \oplus \widehat{\Z_{(p)}}$.  Otherwise, if $\alpha \in \widehat{\Z_{(p)}}$ is a ring element that is not algebraic over $0$, multiplication of the second coordinate by $\alpha$ induces a non-unipotent automorphism $\varphi$ of $K / K^\circ$, and if $\overline{G}$ is a saturated model extending $K$, $\varphi$ can be extended to a non-unipotent automorphism of $\overline{G} / \overline{G}^\circ$.  So by Theorem~\ref{cases}, there is a model of $\Th(G, +)$ of the form $C \oplus D$ where $C$ is a direct sum of finite cyclic groups and $D$ is divisible.

Again, we break into cases according to the form of $C$.  If there is a prime $p$ such that for infinitely many $k \in \omega$, $\Z / {p^k \Z}$ is a direct summand of $C$, then, as in Case~A of 1 $\Rightarrow$ 3 above, we get a contradiction to superstability.  If there are infinitely many distinct primes $p$ for which there is a $k \in \omega$ such that $\left( \Z / {p^k \Z}\right)^{\omega}$ is a direct summand of $C$, again we contradict superstability.  On the other hand, if there are infinitely many distinct primes $p$ for which there is a $k$ such that $\Z / {p^k \Z}$ occurs finitely often as a summand of $C$, then it is simple to produce a non-unipotent automorphism of $C / C^{\circ}$ by combining order-$(p-1)$ permutations of all such summands.  The only case left is that $C$ has bounded exponent.  As noted in the proof of 1 $\Rightarrow$ 3, this implies that $G$ itself is a direct sum of a divisible group and a group of bounded exponent.

\bibliography{SB}

\begin{thebibliography}{10}

\bibitem{baldwin}
John Baldwin.
\newblock {\em Fundamentals of Stability Theory}.
\newblock Perspectives in Mathematical Logic. Springer-Verlag, 1988.

\bibitem{my_thesis}
John Goodrick.
\newblock {\em When are elementarily bi-embeddable models isomorphic?}
\newblock PhD thesis, University of California, Berkeley, 2007.

\bibitem{when}
John Goodrick.
\newblock When does elementary bi-embeddability imply isomorphism? (preprint).
\newblock 2007.

\bibitem{hodges}
Wilfrid Hodges.
\newblock {\em A Shorter Model Theory}.
\newblock Cambridge University Press, 1997.

\bibitem{kaplansky}
Irving Kaplansky.
\newblock {\em Infinite Abelian Groups}.
\newblock The University of Michigan Press, 1954.

\bibitem{nur2}
T.~A. Nurmagambetov.
\newblock The mutual embeddability of models.
\newblock In {\em Theory of Algebraic Structures (in Russian)}, pages 109--115.
  Karagand. Gos. Univ., 1985.

\bibitem{nur1}
T.~A. Nurmagambetov.
\newblock Characterization of $\omega$-stable theories with a bounded number of
  dimensions.
\newblock {\em Algebra i Logika}, 28(5):388--396, 1989.

\bibitem{pillay}
Anand Pillay.
\newblock {\em Geometric Stability Theory}.
\newblock Oxford University Press, 1996.

\bibitem{prest}
Mike Prest.
\newblock {\em Model Theory and Modules}.
\newblock Cambridge University Press, 1988.

\bibitem{bible}
Saharon Shelah.
\newblock {\em Classification Theory}, volume~92 of {\em Studies in logic and
  the foundations of mathematics}.
\newblock North-Holland, second edition, 1990.

\end{thebibliography}

\end{document}